\documentclass[11pt]{amsart}
\usepackage{amscd, amsmath, amsthm, amssymb}

\newcommand{\calF}{\mathcal{F}}

\newcommand{\calO}{\mathcal{O}}

\newcommand{\calR}{\mathcal{R}}

\newcommand{\bbC}{\mathbb{C}}

\newcommand{\bbL}{\mathbb{L}}

\newcommand{\bbP}{\mathbb{P}}

\newcommand{\bbZ}{\mathbb{Z}}

\newcommand\Pic{{\text{Pic}}}

\DeclareMathOperator{\disc}{disc}

\newtheorem{theorem}{Theorem}[section]
\newtheorem{lemma}[theorem]{Lemma}
\newtheorem{proposition}[theorem]{Proposition}

\theoremstyle{definition}     

\theoremstyle{remark}

\numberwithin{equation}{section}

\DeclareMathOperator{\rank}{rank}
\begin{document}
\title[Fake projective plane]
{A Fake Projective Plane with an order 7 automorphism}

\author[J. Keum]{JongHae Keum }
\address{School of Mathematics, Korea Institute for Advanced
Study, Seoul 130-722, Korea } \email{jhkeum@kias.re.kr}
\thanks{Research supported by KOSEF grant R01-2003-000-11634-0}
\begin{abstract}
A fake projective plane is a compact complex surface (a compact complex
manifold of dimension 2) with the same Betti numbers as the complex
projective plane, but not isomorphic to the complex
projective plane.
As was shown by Mumford, there exists at least one such surface.

In this paper we prove the existence of a fake projective plane
which is birational to a cyclic cover of degree 7 of a Dolgachev
surface.
\end{abstract}
\maketitle
\section{Introduction}
It is known that a compact complex surface with the same Betti numbers as the
complex projective plane $\bbP^2$ is projective (see e.g. \cite{BHPV}).
Such a surface is called {\it a fake projective plane} if it is not
isomorphic to $\bbP^2$.

Let $S$ be a fake projective plane, i.e. a surface with
$b_1(S)=0$, $b_2(S)=1$ and $S\ncong \bbP^2$. Then its canonical
bundle is ample and $S$ is of general type. So a fake projective
plane is nothing but a surface of general type with $p_g=0$ and
$c_1(S)^2=3c_2(S)=9$. Furthermore, its fundamental group
$\pi_1(S)$ is infinite. Indeed, by Castelnuovo's rationality
criterion, its second plurigenus $P_2(S)$ must be positive and
hence the first Chern class $-c_1(S)$ of the cotangent bundle of
$S$ can be represented by a K\"ahler form. Then it follows from
the solution of S.-T. Yau \cite{Yau} to Calabi conjecture that $S$
admits a K\"ahler-Einstein metric, and hence its universal cover
is the unit ball in $\bbC^2$.

Mumford \cite{Mum} first discovered a fake projective
 plane. His construction uses the theory of the $p$-adic unit ball by Kurihara \cite{Ku} and
Mustafin \cite {Mus}. Later, using the same idea, Ishida and Kato
\cite{IshidaKato} proved the existence of at least two more. It is not known
if any of these surfaces admits an order 7 automorphism.

In this paper we prove the existence of a fake projective plane with
an order 7 automorphism (see Theorem
\ref{main}). Our construction uses Ishida's
description \cite{Ishida} of an elliptic surface $Y$ covered by a (blow-up) of Mumford's
surface $M$. Recall that there exist an unramified Galois cover $V\to
M$ of degree 8, and a simple group $G$ of order 168 acting on $V$ such
that $Y$ is the minimal resolution of the quotient $V/G$, and $M$
is the quotient of $V$ by a 2-Sylow subgroup of $G$.
Our surface is birationally isomorphic to a cyclic
cover of degree 7 of a cyclic
cover of degree 3 of $Y$.
Mumford's surface $M$ is a degree 21 non-Galois cover of $Y$
\cite{Ishida}, but it is not clear whether it is different from our surface.

The elliptic surface $Y$ has not up to now been constructed directly
(although its properties are stated explicitly), so it does not yet yield
an alternative approach to Mumford's surface. Thus we
do not know how to construct our surface in a direct way.

\medskip
{\it Acknowledgements}

\medskip
I like to thank I. Dolgachev for many useful
conversations. I am also grateful to F. Kato for an e-mail exchange.

\section{Two Dolgachev surfaces}

In \cite{Ishida} Ishida discusses an elliptic surface $Y$ with
$p_g=q=0$ having two multiple fibres of multiplicity 2 and 3
respectively, and proves that the Mumford fake plane $M$ is a
cover of $Y$ of degree 21, but not a Galois cover.

The surface $Y$ is a Dolgachev surface \cite{BHPV}. In particular, it
is simply connected and of Kodaira dimension 1. Besides the two multiple
fibres, its elliptic fibration $|F_Y|$ has
4 more singular fibres $F_1, F_2, F_3, F_4$, all of type $I_3$.
It has also a sextuple section $E$ which is a smooth rational curve meeting one
component of each of $F_1, F_2, F_3$  in 6 points, and two components
of $F_4$ in 1 point and 5 points each.
Write
$$F_i=A_{i1}+A_{i2}+A_{i3}, \quad i=1,2,3,4.$$
After suitable renumbering, one may assume that
$$E\cdot F_Y=E\cdot A_{13}=E\cdot A_{23}=E\cdot A_{33}=6,\quad E\cdot
A_{41}=1,\quad E\cdot A_{43}=5.$$
Note that $E^2=-3$.
The six curves
$$A_{11},\, A_{12},\, A_{21},\, A_{22},\, A_{31},\, A_{32}$$
form a Dynkin diagram of type $A_2^{\oplus 3}$,
$$(-2)\, \textrm{---}\,(-2)\quad (-2)\, \textrm{---}\,(-2)\quad (-2)\,
\textrm{---}\,(-2)$$ and hence can be contracted to 3 singular
points of type $\frac13(1,2)$. Let
$$\sigma\colon Y\to Y'$$
be the contraction morphism. Since $Y$ is simply connected,
$H^2(Y, \bbZ)$ has no torsion and is a lattice under intersection
pairing. Let $$R\subset H^2(Y, \bbZ)$$ be the sublattice generated
by the classes of the six exceptional curves,
 and let $\overline{R}$ and $R^{\perp}$ be its
primitive closure and its orthogonal complement in $H^2(Y, \bbZ)$, respectively.
Since $H^2(Y, \bbZ)$ is unimodular of rank 10, we see that
$$\rank R^{\perp}=4, \quad \disc(\overline{R})\cong -\disc(R^{\perp}).$$
Note that the classes of the 3 curves
$$E,\quad A_{41},\quad A_{42}$$
generates a rank 3 sublattice of $R^{\perp}$ whose discriminant
group has length 1 (a cyclic group of order 7), hence
$$l(\disc(\overline{R}))=l(\disc(R^{\perp}))\le 2.$$
Since $\disc(R)$ is 3-elementary of length $3$, we see that
$\overline{R}$ is an overlattice of index 3
of $R$. This implies that there is a cyclic cover of degree 3
$$X\to Y'$$
branched exactly at the 3 singular points of $Y'$. Then $X$ is a
nonsingular surface. It turns out that $X$ is another Dolgachev
surface. More precisely, we have

\begin{proposition}\label{x} The surface $X$ is a minimal elliptic
surface of Kodaira dimension $1$
with $p_g=q=0$ with one fibre of multiplicity $2$, one of
multiplicity $3$,
one singular fibre of type $I_9$ and three of type $I_1$.
Furthermore, it has $3$ sextuple sections $E_1, E_2, E_3$ which
together with $6$ components of the fibre of type $I_9$ can be contracted to
$3$ singular points of type $\frac17(1,3)$.
\end{proposition}

\begin{proof} The elliptic fibration $|F_Y|$ on $Y$ induces an elliptic
fibration of $Y'$, as the six exceptional curves are contained in fibres. The
image $\sigma(E)$ of the $-3$-curve $E$ is a smooth rational curve
not passing through any of the 3 singular points of $Y'$, hence splits
in $X$ to give 3 smooth rational curves $E_1, E_2, E_3$. This implies
that the fibres of $Y'$ do not split in $X$. The fibre containing one
of the singular points of $Y'$ gives a fibre of type $I_1$, the fibre of
type $I_3$ gives a fibre of type $I_9$, and the multiple fibres give
multiple fibres of the same multiplicities. Adding up the Euler number
of fibres, we have $c_2(X)=12$. Thus $X$ is minimal. This proves the first
assertion.

Note that the 3 curves $\sigma(E), \sigma(A_{41}), \sigma(A_{42})$ on $Y'$ form
a configuration of smooth rational curves
$$(-3)\,\textrm{---}\,(-2)\,\textrm{---}\,(-2)$$
which may be contracted to a singular point of type $\frac17(1,3)$.
Clearly their preimages in $X$ form a disjoint union of 3 such configurations.
\end{proof}

\section{Construction of a fake projective plane}

Let $X$ be the surface from Proposition \ref{x}. Denote by $|F|$ the
elliptic fibration on $X$ induced from $|F_Y|$ on $Y$.
Denote clockwise the components of the fibre of $|F|$ of type $I_9$  by
$$A_1,\, A_2,\, A_3,\, B_1,\, B_2,\, B_3,\, C_1,\, C_2,\, C_3$$
in such a way that
$$E_1\cdot A_2= E_2\cdot B_2= E_3\cdot C_2=1.$$
Thus the 9 curves
$$A_1, A_2, E_1,\quad B_1, B_2, E_2,\quad C_1, C_2, E_3$$
form a dual diagram
$$(-2)\,\textrm{---}\,(-2)\,\textrm{---}\,(-3)\quad (-2)\,\textrm{---}\,(-2)\,\textrm{---}\,(-3)\quad (-2)\,\textrm{---}\,(-2)\,\textrm{---}\,(-3)$$
which can be contracted to 3 singular points of type $\frac17(1,3)$.
Let
$$\nu\colon X\to X'$$
be the birational contraction morphism. Let $$\calR\subset H^2(X,
\bbZ)\cong \Pic(X)$$  be the sublattice generated by the classes
of the nine exceptional curves,
 and let $\overline{\calR}$ and $\calR^{\perp}$ be its
primitive closure and its orthogonal complement in $H^2(X, \bbZ)$, respectively.
Since $H^2(X, \bbZ)$ is unimodular of rank 10, we see that
$\rank \calR^{\perp}=1$ and hence
$$l(\disc(\overline{\calR}))=l(\disc(\calR^{\perp}))=1.$$
Since $\disc(\calR)$ is 7-elementary of length $3$, we see that
$\overline{\calR}$ is an overlattice of index 7 of $\calR$. This
implies that there is a cyclic cover of degree 7
$$\pi\colon Z\to X'$$
branched exactly at the 3 singular points of $X'$.
Then $Z$ is a nonsingular surface.

\begin{theorem} \label{main} The surface $Z$ is a fake projective
plane, i.e. a surface of general type with $p_g(Z)=0$ and $c_2(Z)=3$.
\end{theorem}

First we show that $Z$ is of general type.

\begin{lemma} \label{gt} The surface $Z$ is of general type.
\end{lemma}

\begin{proof}
Let $X^0$ be the smooth part of $X'$, and $\kappa(X^0)$
be the Kodaira dimension of $X^0$ (that is, its logarithmic
Kodaira dimension). Then
$$\kappa(Z)\ge \kappa(X^0)\ge \kappa(X)=1.$$
We know that $X'$ has Picard number 1,
  hence is relatively minimal, i.e. contains no curve $C$ with
$C\cdot K_{X'}<0$ and $C^2<0$.

Suppose $\kappa(X^0)=1$. Then there is an elliptic fibration on $X'$
(\cite{Ka} Theorem 2.3, \cite{Mi} Ch.II, Theorem 6.1.4, \cite{KZ}
Theorem 4.1). This implies that $X$ admits an elliptic fibration whose
fibres contain the nine exceptional curves, but no elliptic fibration
on $X$ can contain a $-3$-curve in its fibres because $X$ is
minimal. Thus $\kappa(X^0)=2$ and the
assertion follows.
\end{proof}

\begin{lemma} \label{c2} $c_2(Z)=3$.
\end{lemma}

\begin{proof} Let $e(X^0)$ be the Euler number of $X^0$. Then
$$c_2(Z)=7e(X^0)+3=7\{e(X)-12\}+3=3.$$
\end{proof}

Let $$E_1\cdot A_3=\alpha,\quad E_1\cdot B_3=\beta,\quad E_1\cdot
C_3=\gamma.$$ Then
\begin{equation}\label{1}
\alpha +\beta +\gamma=5.
\end{equation}
We claim that there are two possible cases for
$(\alpha,\beta,\gamma)$;

\medskip
Case I: $(\alpha,\beta,\gamma)=(2,1,2)$.

Case II: $(\alpha,\beta,\gamma)=(1,3,1)$.

\medskip
Let us prove the claim. Since $X$ admits an order 3 automorphism
which rotates the $I_9$ fibre and the 6-sections $E_1, E_2, E_3$,
we have
$$E_2\cdot (A_3, B_3,C_3)=(\gamma, \alpha, \beta),\quad E_3\cdot (A_3, B_3,
C_3)=(\beta, \gamma, \alpha).$$ The divisor class $E_2-E_1$ is
orthogonal to the class $F$ of a fibre of the elliptic fibration
on $X$ and hence can be written in the form
$$E_2-E_1=\sum_{i=1}^{3}a_iA_i +\sum_{i=1}^{3}b_iB_i +\sum_{i=1}^{3}c_iC_i$$
for some rational numbers $a_i, b_i, c_i$. Applying the order 3
automorphism, we get
$$E_3-E_2=\sum_{i=1}^{3}a_iB_i +\sum_{i=1}^{3}b_iC_i +\sum_{i=1}^{3}c_iA_i,$$
$$E_1-E_3=\sum_{i=1}^{3}a_iC_i +\sum_{i=1}^{3}b_iA_i +\sum_{i=1}^{3}c_iB_i.$$
Adding the three equations side by side, we get
\begin{equation}\label{2}
a_i+b_i+c_i=0,\quad i=1,2,3.
\end{equation}
Intersecting $E_2-E_1$ with $A_i, B_i, C_i, E_i$, we get 12
equations in 9 unknowns $a_i, b_i, c_i$. The system of these 12
equations together with the 3 equations from \eqref{2} has a
solution if and only if $(\alpha,\beta,\gamma)=(2,1,2)$ or
$(1,3,1)$. This completes the proof of the claim.

By a direct calculation, we see that the discriminant group of $\calR$
is generated by the three elements
$$\frac17(A_1+2A_2+3E_1),\quad \frac17(B_1+2B_2+3E_2),\quad \frac17(C_1+2C_2+3E_3)$$
and hence its discriminant form is isomorphic to $(-3/7)^{\oplus 3}$.
Thus a generator of the quotient group  $\overline{\calR}/\calR$ is
of the form
$$v=\frac17(A_1+2A_2+3E_1)+\frac{a}7(B_1+2B_2+3E_2)+\frac{b}7(C_1+2C_2+3E_3).$$
Since the intersection numbers $v\cdot A_3$ and $v\cdot B_3$ are
integers, we have $a\equiv 4$ (mod 7) and $b\equiv 2$ (mod 7) in
Case I, and $a\equiv 2$ (mod 7) and $b\equiv 4$ (mod 7) in Case
II. This determines $v$ uniquely modulo $\calR$ in each case. We
fix an effective divisor
$$B=\left\{\begin{array}{l}A_1+2A_2+3E_1+4B_1+B_2+5E_2+2C_1+4C_2+6E_3\,\,({\rm Case}\, {\rm I}) \\
A_1+2A_2+3E_1+2B_1+4B_2+6E_2+4C_1+C_2+5E_3\,\,({\rm Case}\, {\rm
II})\end{array}\right.$$ Then $B$ is divisible by 7 in $\Pic(X)$,
so we can write
$$\calO_X(B)\cong \calO_X(7L)$$ for some divisor $L$ on $X$.
Since $X$ is simply connected, $\Pic(X)$ has no torsion and $L$ is
determined uniquely by $B$ up to linear equivalence.
We denote by $\bbL$ the total space of $\calO_X(L)$ and by
$$p\colon\bbL\to X$$ the bundle projection. If $t\in H^0(\bbL, p^*\calO_X(L))$ is the
tautological section, and if $s\in H^0(X, \calO_X(B))$ is the
section vanishing exactly along $B$, then the zero divisor of
$p^*s-t^7$ defines an analytic subspace $W$ in $\bbL$,
$$W:=(p^*s-t^7=0)\subset \bbL.$$ Since $B$ is not reduced, $W$ is not normal.

\begin{lemma} \label{w} $H^2(W, \calO_W)=0$.
\end{lemma}

\begin{proof} Let $$p\colon W\to X$$ be the restriction to $W$ of the bundle
projection $p\colon \bbL\to X$. Since it is a finite morphism,
$$H^2(W, \calO_W)=H^2(X, p_*\calO_W).$$
We know that
$$p_*\calO_W=\calO_X\oplus \calO_X(-L)\oplus \calO_X(-2L)
\oplus ...\oplus \calO_X(-6L).$$
Thus
$$H^2(W, \calO_W)=\bigoplus_{i=0}^6 H^0(X,\calO_X(K_X+iL)).$$
We know $H^0(X,\calO_X(K_X))=0.$

Assume $K_X+iL$ $(1\le i\le 6)$ is effective.

\subsection*{Case I: $E_1\cdot (A_3, B_3,C_3)=(2,1,2)$.}

\medskip\noindent

To get a contradiction we use the following intersection
numbers:
$$\renewcommand{\arraystretch}{1.3}\begin{array}{l}
L\cdot A_1=L\cdot A_2=0, \qquad L\cdot A_3 =4,\\
L\cdot B_1=-1, \quad L\cdot B_2=1, \quad L\cdot B_3 =4,\\
L\cdot C_1=L\cdot C_2=0, \qquad L\cdot C_3 =4,\\
L\cdot E_1=-1, \qquad L\cdot E_2=L\cdot E_3 =-2,\\
L\cdot K_X=2, \qquad L\cdot F=12.
\end{array}$$
We claim that there are nonnegative integers $i_1, i_2, i_3$ with
$i_1+i_2+i_3=2i+1$  and an
effective divisor $G_i$ with support contained in the fibre of
type $I_9$ such that the divisor
$$D_i=K_X+iL-i_1E_1-i_2E_2-i_3E_3-G_i$$
is effective. This contradicts the fact that
$F$ is represented by an irreducible curve with self-intersection 0,
as we have
$$D_i\cdot F= iL\cdot F-6(i_1+i_2+i_3)=12i-6(2i+1)=-6<0.$$
This proves that $H^0(X,\calO_X(K_X+iL))=0$.

It remains to prove the claim.

\medskip
Assume $i=6$.

\medskip\noindent

The divisor $D_1'=K_X+6L- E_1-E_2-E_3$ is effective, because
$$(K_X+6L)\cdot E_1=-5,\quad (K_X+6L)\cdot E_2=-11,\quad (K_X+6L)\cdot
E_3=-11.$$ Since $D_1'\cdot E_3<0$, $D_1'-E_3=K_X+6L-
E_1-E_2-2E_3$ is effective. Iterating this process, we see that
the divisor
$$D_6=K_X+6L-3E_1-5E_2-5E_3-(A_1+2A_2+A_3+4B_1+B_2+C_1+3C_2)$$
is effective.

\medskip\noindent
The other cases $i=5,4,3,2,1$ can be handled similarly. We give
$D_i$ in each case for the sake of completeness:
$$\begin{array}{lr}
D_5=K_X+5L-3E_1-4E_2-4E_3 \\
 & \kern-4.5cm -(A_1+2A_2+A_3+3B_1+B_2+C_1+2C_2),\\[4pt]
D_4=K_X+4L-2E_1-3E_2-4E_3 \\
 & \kern-4.5cm -(A_1+A_2+2B_1+B_2+C_1+2C_2),\\[4pt]
D_3=K_X+3L-E_1-3E_2-3E_3 \\
 & \kern-4.5cm -(A_1+A_2+2B_1+B_2+B_3+C_1+2C_2),\\[4pt]
D_2=K_X+2L-E_1-E_2-3E_3 \\
 & \kern-4.5cm -(A_1+A_2+B_1+C_1+2C_2+C_3),\\[4pt]
D_1=K_X+L-E_1-E_2-E_3 \\
 & \kern-4.5cm -(B_1+B_2+B_3+C_1+C_2).
\end{array}$$

\subsection*{Case II: $E_1\cdot (A_3, B_3,C_3)=(1,3,1)$.}

\medskip\noindent

In this case we use the following intersection numbers:
$$\renewcommand{\arraystretch}{1.3}\begin{array}{l}
L\cdot A_1=L\cdot A_2=0, \qquad L\cdot A_3 =4,\\
L\cdot B_1=L\cdot B_2=0, \quad L\cdot B_3 =4,\\
L\cdot C_1=-1, \quad L\cdot C_2=1, \qquad L\cdot C_3 =4,\\
L\cdot E_1=-1, \qquad L\cdot E_2=L\cdot E_3 =-2,\\
L\cdot K_X=2, \qquad L\cdot F=12.
\end{array}$$
The rest of the proof is similar to the previous case.
\end{proof}

Let
$$\mu\colon W'\to W$$
be the normalization, and $$f\colon \widetilde{W}\to W'$$ the
minimal resolution. Then $\widetilde{W}$ is a nonsingular surface
birational to $Z$. In particular, $$p_g(Z)=p_g(\widetilde{W}).$$

\begin{lemma} \label{pg} $p_g(Z)=0$.
\end{lemma}

\begin{proof}
Consider the exact sequence of sheaves on $W$
$$0\to \calO_W\to \mu_*\calO_{W'} \to \calF \to 0.$$
Since $\calF$ is supported in dimension 1, $$H^2(W, \calF)=0.$$ By
Lemma \ref{w}, $$H^2(W,\calO_{W})=0.$$ Thus the long exact
sequence gives $$H^2(W,\mu_*\calO_{W'})=0.$$ Since normalization
is an affine morphism, this implies that
 $$H^2(W',\calO_{W'})=H^2(W,\mu_*\calO_{W'})=0.$$
Since $W'$ is normal, $$f_*\calO_{\widetilde{W}}=\calO_{W'},$$ and
hence
$$H^2(W',f_*\calO_{\widetilde{W}})=H^2(W',\calO_{W'})=0.$$
Note that the direct image sheaf $R^1f_*\calO_{\widetilde{W}}$ is
supported in dimension 0. Thus
$$H^1(W',R^1f_*\calO_{\widetilde{W}})=0.$$ Now the Leray spectral
sequence gives
$$p_g(\widetilde{W})=H^2(\widetilde{W},\calO_{\widetilde{W}})=0.$$
\end{proof}

Now Theorem \ref{main} follows from Lemma \ref{gt},  \ref{c2} and \ref{pg}.


\end{document}